\documentclass[12pt,thmsa,sw20bams]{article}

\begin{document}

\begin{titlepage}
\title{\bf Lifts of Almost r-Contact and r-Paracontact Structures}
\author{ Mehmet Tekkoyun \footnote{tekkoyun@pau.edu.tr} \\
 {\small Department of Mathematics, Pamukkale University,}\\
{\small 20070 Denizli, Turkey}}
\date{\today}
\maketitle
\begin{abstract}

In this study, taking into considering lifting theory, we shall
obtain both almost complex and paracomplex structures on the
tangent bundle, based on almost Lorentzian r-contact and
r-paracontact manifold.


 {\bf
M.S.C.2000:} 53D15, 53B30, 28A51.
\end{abstract}

\section{\textbf{Introduction}}

In modern differentiable geometry, lifting theory has a considerable
position. Because, using lift method it is possible to generalize to
differentiable structures on any manifold to its extensions. Also, Vertical,
complete and horizontal lifts of functions, vector fields, 1-forms and other
tensor fields defined on any (complex) manifold $M$ to tangent bundle $TM$
has been obtained and the basic classical results can be traced in the
well-known papers [1-8]. It is introduced almost r-contact structures on the
tangent bundle [9].

The paper is structured as follows. In section 2, we recall complex,
paracomplex, contact and Lorentzian differential structures [9-14]. In
section 3, we give vertical, complete and horizontal lifts of complex (or
paracomplex) structures [9]. In section 4, we produce almost complex (or
paracomplex) structures by means of complete and horizontal lifts of almost
Lorentzian r-contact (or $r-$paracontact) structure on almost Lorentzian
r-contact ( or $r-$paracontact) manifold. In the finally section, we see as
to what kind of structure is defined on the tangent bundle $T(M)$ when we
consider an almost Lorentzian $r-$ contact (or $r-$paracontact) structure on
the base manifold.

Along this paper, all mappings and manifolds will be understood to be of
class differentiable and the sum is taken over repeated indices. In this
study we denote by $\Im _{s}^{r}(M)$ the set of all tensor fields of class $%
C^{\infty }$ and of type $(r,s)$ in $M$. We now put $\Im
(M)=\sum_{r,s=1}^{\infty }\Im _{s}^{r}(M)$; which is the set of all tensor
fields in $M$. For example $\Im _{0}^{0}(M),\Im _{0}^{1}(M),\Im _{1}^{0}(M)$
and $\Im _{1}^{1}(M)$ are the set of functions, vector fields,1-forms and
tensor fields of type (1,1) on $M$, respectively. We give by $\mathcal{F}(M)$
to $\Im _{0}^{0}(M),$ by $\chi (M)$ to $\Im _{0}^{1}(M)$ and by $\chi
^{*}(M) $ to $\Im _{1}^{0}(M).$ Similarly, we respectively denote by $\Im
_{s}^{r}(T(M))$ and $\Im (T(M))$ the corresponding sets of tensor fields in
the tangent bundle $T(M)$. Also $v$, $c$ and $h$ will denote the vertical,
complete and horizontal lifts to $TM$ of geometric structures on $M$,
respectively.

\subsection{Complex and Paracomplex manifolds}

Let $TM$ be tangent bundle of $m$- dimensional a differential manifold $\,M.$
A tensor field $J$ on $TM$ is called an \textit{almost complex (or
paracomplex) structure} on $TM$ if at every point $p$ of $TM,$ $J$ is
endomorphism of the tangent space $T_{p}M$ such that $J^{2}=-I$ (or $J^{2}=I$%
), i.e., shortly it can shown by $J^{2}=\epsilon I.$ For $\epsilon =-1$ (or $%
\epsilon =1$)$,$ a manifold $M$ with fixed almost complex (or paracomplex)
structure $J$ is called \textit{almost complex (or paracomplex) manifold}.

Let $(x^{i},\,y^{i})$ be a real coordinate system on a neighborhood $U$ of
any point $p$ of $TM$. In this case, it is respectively defined by $\{\frac
{\partial }{\partial x^{i}}|_{p},\frac {\partial }{\partial y^{i}}|_{p}\}$
and $\{\left. dx^{i}\right| _{p},\left. dy^{i}\right| _{p}\}$ natural bases
over $\mathbf{R}$ of tangent space $T_{p}(M)$ and cotangent space $%
T_{p}^{*}(M)$ of $M.$ Then $TM$ is called \textit{\ complex ( or
paracomplex) manifold} if there exists an open covering $\left\{ U\right\} $
of $M$ satisfying the following condition: There is a local coordinate
system $(x^{i},\,y^{i})$ on each $U,$ such that

\begin{equation}
J(\frac{\partial }{\partial x^{i}})=\frac{\partial }{\partial y^{i}},\,J(%
\frac{\partial }{\partial y^{i}})=\epsilon \frac{\partial }{\partial x^{i}},
\label{1.1}
\end{equation}

for each point of $U.$

Let $(z^{i},\overline{z}^{i})$ be an complex (or paracomplex) local
coordinate system on a neighborhood $U$ of any point $p$ of $M$ where $%
z^{i}=x^{i}+${\textbf{i}$\,y^{i},\,$\textbf{i}$^{2}=\epsilon.$ Therefore
this coordinates are defined by locally $x^{i},\,y^{i}:M\rightarrow \mathbf{R%
}^{m},\forall A,\overline{A}\in M,$ $A=a^{i}+$\textbf{i}$b^{i},$ $%
x^{i}(A)=a^{i}$, $y^{i}(A)=b^{i}$ and $z^{i},\overline{z}^{i}:M\rightarrow
\mathbf{C}_{\epsilon }^{m},$ $z^{i}(A)=A$, $\overline{z}^{i}(A)=\overline{A}$%
, where $\mathbf{R}$ is real number, $\mathbf{C}_{\epsilon }$ is complex (or
paracomplex) numbers$.$ }

We define the vector fields
\begin{equation}
\frac {\partial }{\partial z^{i}}|_{p}=\frac{1}{2}\{\frac {\partial
}{\partial x^{i}}|_{p}-\mathbf{i}\frac {\partial }{\partial
y^{i}}|_{p}\},\frac {\partial }{\partial \overline{z}^{i}}|_{p}=\frac{1}{2}\{%
\frac{\partial }{\partial x^{i}}|_{p}+\mathbf{i}\frac {\partial }{\partial
y^{i}}|_{p}\}  \label{1.2}
\end{equation}

and the dual covector fields

\begin{equation}
\left. dz^{i}\right| _{p}=\left. dx^{i}\right| _{p}+\mathbf{i}\left.
dy^{i}\right| _{p},\left. d\overline{z}^{i}\right| _{p}=\left. dx^{i}\right|
_{p}-\mathbf{i}\left. dy^{i}\right| _{p}  \label{1.3}
\end{equation}

which represent bases of the tangent space $T_{p}(M)$ and cotangent space $%
T_{p}^{*}(M)$ of $M$ respectively. Then the endomorphism $J$ is given as:
\begin{equation}
J(\frac{\partial }{\partial z^{i}})=-\mathbf{i}\epsilon \frac{\partial }{%
\partial z^{i}},J(\frac{\partial }{\partial \overline{z}^{i}})=\mathbf{i}%
\epsilon \frac{\partial }{\partial \overline{z}^{i}}  \label{1.4}
\end{equation}

The dual endomorphism $J^{*}$ of the cotangent space $T_{p}^{*}(M)$ of
manifold $M$ satisfies $J^{*2}=\epsilon I$ and is defined by
\begin{equation}
J^{*}(dz^{i})=-\epsilon \mathbf{i}dz^{i},\,J^{*}(d\overline{z}^{i})=\epsilon
\mathbf{i}d\overline{z}^{i}.  \label{1.5}
\end{equation}

\subsection{Almost Contact and Paracontact Manifolds}

Let $\overline{M}$ be an $n$-dimensional differentiable manifold. If there
exist on $\overline{M}$ a (1,1)-tensor field $\varphi $, a vector field $\xi
$ and a 1-form $\eta $ satisfying

\begin{equation}
\eta (\xi )=1,\,\,\,\varphi ^{2}=\epsilon I-\eta \otimes \xi ,  \label{1.6}
\end{equation}

where $I$ is the identity, then $\overline{M}$ is said to be an almost
contact (or paracontact) manifold for $\epsilon =-1$ (or $\epsilon =1$)$.$
In the almost contact (or paracontact) manifold, the following relations
hold good:

\begin{equation}
\,\varphi (\xi )=0,\,\eta \circ \varphi =0,\,\,rank(\varphi )=n-1.
\label{1.7}
\end{equation}

Every almost contact ( or paracontact) manifold has a positive definite
Riemannian metric $G$ such that

\begin{equation}
\eta (\overline{X})=G(\xi ,\overline{X}),\,\,G(\varphi \overline{X},\varphi
\overline{Y})=G(\overline{X},\overline{Y})-\eta (\overline{X})\eta (%
\overline{Y}),\,\,\overline{X},\overline{Y}\in \chi (\overline{M})
\label{1.8}
\end{equation}

where $\chi (\overline{M})$ denotes the set of differentiable vector fields
on $\overline{M}$. In this case, we say that $\overline{M}$ has an almost
contact (or paracontact) Riemannian structure $(\varphi ,\xi ,\eta ,G)$ and $%
\overline{M}$ is said to be an almost contact (or paracontact) Riemannian
manifold.

This structure generalizes as follows.

\subsection{Almost r-Contact and r-Paracontact Manifolds}

One takes care of a tensor field $\mathbf{F}$ of type (1,1) on a manifold $%
\overline{M}$ of dimension $(2n+r)$. If there exists on $\overline{M}$ the
vector fields $(\xi _{\alpha })$ and the 1- forms $(\eta ^{\alpha })$ such
that

\begin{equation}
\eta ^{\alpha }(\xi _{\beta })=\delta _{\beta }^{\alpha },\,\,\mathbf{F}(\xi
_{\alpha })=0,\,\eta ^{\alpha }\circ {\mathbf{F}}=0,\,\,\,{\mathbf{F}}%
^{2}=\epsilon I+\sum_{\alpha =1}^{r}\xi _{\alpha }\otimes \eta ^{\alpha },
\label{1.9}
\end{equation}

then the structure $(\mathbf{F},\xi _{\alpha },\eta ^{\alpha })$ is an
almost $r-$contact (or $r-$paracontact) structure, where $(\alpha ,\beta
=1,2,...,r)$ and $\delta _{\beta }^{\alpha }$ denotes Kronecker delta. A
manifold $\overline{M}$ endowed with $(\mathbf{F},\xi _{\alpha },\eta
^{\alpha })$-structure is called an almost $r$- contact (or $r$-
paracontact) manifold.

\subsection{ Almost Lorentzian Contact and  Paracontact Manifolds}

Let $\overline{M}$ be an $n$-dimensional differentiable manifold equipped
with a triple $(\varphi ,\xi ,\eta ),$ where $\varphi $ is a (1,1)-tensor
field$,\xi $ a vector field and $\eta $ is a 1-form on $\overline{M}$
satisfying

\begin{equation}
\eta (\xi )=-1,\,\,\,\varphi ^{2}=\epsilon I+\eta \otimes \xi ,  \label{1.10}
\end{equation}

where $I$ is the identity. In manifold $\overline{M},$ the following
relations hold good:

\begin{equation}
\,\varphi (\xi )=0,\,\eta \circ \varphi =0,\,\,rank(\varphi )=n-1.
\label{1.11}
\end{equation}

If $\overline{M}$ has a Lorentzian metric $G$ such that

\begin{equation}
G(\varphi \overline{X},\varphi \overline{Y})=G(\overline{X},\overline{Y}%
)+\eta (\overline{X})\eta (\overline{Y}),\,\,\overline{X},\overline{Y}\in
\chi (\overline{M})  \label{1.12}
\end{equation}

we say that $\overline{M}$ has a Lorentzian almost contact (or paracontact)
structure $(\varphi ,\xi ,\eta ,G)$ and $\overline{M}$ is said to be a
Lorentzian almost contact (or paracontact) manifold.

This structure extends the following as:

\subsection{Almost Lorentzian r-Contact and r-Paracontact Manifolds}

One takes care of a tensor field $\mathbf{F}$ of type (1,1) on a manifold $%
\overline{M}$ of dimension $(2n+r)$. If there exists on $\overline{M}$ the
vector fields $(\xi _{\alpha })$ and the 1- forms $(\eta ^{\alpha })$ such
that

\begin{equation}
\eta ^{\alpha }(\xi _{\beta })=-\delta _{\beta }^{\alpha },\,\,\mathbf{F}%
(\xi _{\alpha })=0,\,\eta ^{\alpha }\circ {\mathbf{F}}=0,\,\,\,{\mathbf{F}}%
^{2}=\epsilon I-\sum_{i=1}^{r}\xi _{\alpha }\otimes \eta ^{\alpha },
\label{1.13}
\end{equation}

then the structure $(\mathbf{F},\xi _{\alpha },\eta ^{\alpha })$ is also
called a Lorentzian almost $r$-contact (or $r-$paracontact) structure, where
$(\alpha ,\beta =1,2,...,r)$ and $\delta _{\beta }^{\alpha }$ denotes
Kronecker delta. A manifold $\overline{M}$ endowed with $(\mathbf{F},\xi
_{\alpha },\eta ^{\alpha })$-structure is called a Lorentzian almost $r$%
-contact (or $r-$paracontact) manifold.

\section{\textbf{Lifts of Almost Lorentzian r-Contact and r-Paracontact
Structures}}

In this section, we obtain almost complex (or paracomplex) structures on
tangent bundle $T(\overline{M})$ of almost Lorentzian $r$- contact (or $r$-
paracontact) manifold $\overline{M}$ having the structure $({\mathbf{F}},\xi
_{\alpha },\eta ^{\alpha }).$

\subsection{Complete Lifts}

\textbf{Theorem 4.1.} Let $\overline{M}$ be a differentiable manifold
endowed with almost $r$- contact (or $r$- paracontact) structure $({\mathbf{F%
}},\xi _{\alpha },\eta ^{\alpha })$, then
\[
\widetilde{J}={\mathbf{F}}^{c}+\sum_{\alpha =1}^{r}\xi _{\alpha }^{v}\otimes
\eta ^{\alpha v}-\xi _{\alpha }^{c}\otimes \eta ^{\alpha c}
\]

is almost complex (or paracomplex) structure on $T(\overline{M})$.

\textbf{Proof :} From (\ref{1.9}) and the vertical and complete lifts of
complex (or paracomplex) tensor fields we have
\begin{equation}
({\mathbf{F}}^{2})^{c}=(\epsilon I+\sum_{\alpha =1}^{r}\xi _{\alpha }\otimes
\eta _{\alpha })^{c},  \label{2.1}
\end{equation}

\begin{equation}
({\mathbf{F}}^{c})^{2}=\epsilon I+\sum_{\alpha =1}^{r}\xi _{\alpha
}^{v}\otimes \eta ^{\alpha c}+\xi _{\alpha }^{c}\otimes \eta ^{\alpha v},
\label{2.2}
\end{equation}

and
\begin{equation}
{\mathbf{F}}^{c}(\xi _{\alpha }^{v})=0,{\mathbf{F}}^{c}(\xi _{\alpha
}^{c})=0,  \label{2.3}
\end{equation}
\begin{equation}
\eta ^{\alpha v}\circ {\mathbf{F}}^{c}=0,\eta ^{\alpha c}\circ {\mathbf{F}}%
^{v}=0,\eta ^{\alpha c}\circ {\mathbf{F}}^{c}=0,  \label{2.4}
\end{equation}
\begin{equation}
\eta ^{\alpha v}(\xi _{\beta }^{v})=0,\eta ^{\alpha v}(\xi _{\beta
}^{c})=\delta _{\beta }^{\alpha },\eta ^{\alpha c}(\xi _{\beta }^{v})=\delta
_{\beta }^{\alpha },\,\eta ^{\alpha c}(\xi _{\beta }^{c})=0.  \label{2.5}
\end{equation}

Consider a structure $\widetilde{J}$ of $\Im _{1}^{1}(T\overline{M})$ given
by

\begin{equation}
\widetilde{J}={\mathbf{F}}^{c}+\sum_{\alpha =1}^{r}(\xi _{\alpha
}^{v}\otimes \eta ^{\alpha v}-\xi _{\alpha }^{c}\otimes \eta ^{\alpha c}).
\label{2.6}
\end{equation}

Using (\ref{2.2}) and (\ref{2.6}) we find equation
\begin{eqnarray}
(\tilde{J})^{2} &=&({\mathbf{F}}^{c}+\sum_{\alpha =1}^{r}\xi _{\alpha
}^{v}\otimes \eta ^{\alpha v}-\xi _{\alpha }^{c}\otimes \eta ^{\alpha c})^{2}
\label{2.7} \\
&=&\epsilon I+\sum_{\alpha =1}^{r}[\xi _{\alpha }^{v}\otimes \eta ^{\alpha
c}+\xi _{\alpha }^{c}\otimes \eta ^{\alpha v}++{\mathbf{F}}^{c}(\xi _{\alpha
}^{v})\eta ^{\alpha v}-{\mathbf{F}}^{c}(\xi _{\alpha }^{c})\eta ^{\alpha c}+
\nonumber \\
&&+(\eta ^{\alpha v}\circ {\mathbf{F}}^{c})\xi _{\alpha }^{v}-(\eta ^{\alpha
c}\circ {\mathbf{F}}^{c})\xi _{\alpha }^{c}+\xi _{\alpha }^{v}\otimes (\eta
^{\alpha v}(\xi _{\alpha }^{v}))\eta ^{\alpha v}-\xi _{\alpha }^{v}\otimes
(\eta ^{\alpha v}(\xi _{\alpha }^{c}))\eta ^{\alpha c}  \nonumber \\
&&-\xi _{\alpha }^{c}\otimes (\eta ^{\alpha c}(\xi _{\alpha }^{v}))\eta
^{\alpha v}+\xi _{\alpha }^{c}\otimes (\eta ^{\alpha c}(\xi _{\alpha
}^{c}))\eta ^{\alpha c}].  \nonumber
\end{eqnarray}
By means of (\ref{2.3}), (\ref{2.4}) and (\ref{2.5}), we have

\begin{equation}
(\widetilde{J})^{2}=\epsilon I.  \label{2.8}
\end{equation}
So, $\widetilde{J}$ is an almost complex (or paracomplex) structure in $T(%
\bar{M})$. Hence the proof is completed.

Thinking lift properties of tensor fields and the equation (\ref{2.6}), we

\[
\widetilde{J}X^{v}=({\mathbf{F}}X)^{v}-(\eta ^{\alpha }(X))^{v}U_{\alpha
}^{c}
\]

\[
\widetilde{J}X^{c}=({\mathbf{F}}X)^{c}+(\eta ^{\alpha }(X))^{v}U_{\alpha
}^{v}-(\eta ^{\alpha }(X))^{c}U_{\alpha }^{c}.
\]

where $X\in \chi (\overline{M}),X^{v},X^{c}\in \chi (T\overline{M}).$

For example, if $\eta ^{\alpha }(X)=0,$ we obtain

\[
\widetilde{J}X^{v}=({\mathbf{F}}X)^{v},\,\widetilde{J}X^{c}=({\mathbf{F}}%
X)^{c}.
\]

In view of the structure $({\mathbf{F}},\xi _{\alpha },\eta ^{\alpha })\,$%
and the equations (\ref{2.3}),(\ref{2.5}) and (\ref{2.6}) it follows

\[
\widetilde{J}\xi _{\alpha }^{v}=-\delta _{\beta }^{\alpha }\xi _{\alpha
}^{c}=\xi _{\beta }^{c}\,,\,\,\,\,\,\widetilde{J}\xi _{\alpha }^{c}=\delta
_{\beta }^{\alpha }\xi _{\alpha }^{v}=\xi _{\beta }^{v}\,,\,\,\alpha ,\beta
=1,2,...,r.
\]

\textbf{Theorem 4.2.} Let $\overline{M}$ be a differentiable manifold
endowed with almost Lorentzian $r-$contact ( or $r-$paracontact) structure $(%
{\mathbf{F}},\xi _{\alpha },\eta ^{\alpha })$, then an almost complex (or
paracomplex) structure on $T(\overline{M})$ is calculated by
\[
\widehat{J}={\mathbf{F}}^{c}-\sum_{\alpha =1}^{r}\xi _{\alpha }^{v}\otimes
\eta ^{\alpha v}-\xi _{\alpha }^{c}\otimes \eta ^{\alpha c}.
\]

\textbf{Proof :} By means of the equation (\ref{1.13}) and the vertical and
complete lifts of complex (or paracomplex) tensor fields we have
\begin{equation}
({\mathbf{F}}^{2})^{c}=(\epsilon I-\sum_{\alpha =1}^{r}\xi _{\alpha }\otimes
\eta _{\alpha })^{c},  \label{2.9}
\end{equation}

\begin{equation}
({\mathbf{F}}^{c})^{2}=\epsilon I-\sum_{\alpha =1}^{r}\xi _{\alpha
}^{v}\otimes \eta ^{\alpha c}+\xi _{\alpha }^{c}\otimes \eta ^{\alpha v},
\label{2.10}
\end{equation}

and
\begin{equation}
{\mathbf{F}}^{c}(\xi _{\alpha }^{v})=0,{\mathbf{F}}^{c}(\xi _{\alpha
}^{c})=0,  \label{2.11}
\end{equation}
\begin{equation}
\eta ^{\alpha v}\circ {\mathbf{F}}^{c}=0,\eta ^{\alpha c}\circ {\mathbf{F}}%
^{v}=0,\eta ^{\alpha c}\circ {\mathbf{F}}^{c}=0,  \label{2.12}
\end{equation}
\begin{equation}
\eta ^{\alpha v}(\xi _{\beta }^{v})=0,\eta ^{\alpha v}(\xi _{\beta
}^{c})=-\delta _{\beta }^{\alpha },\eta ^{\alpha c}(\xi _{\beta
}^{v})=-\delta _{\beta }^{\alpha },\,\eta ^{\alpha c}(\xi _{\beta }^{c})=0.
\label{2.13}
\end{equation}

Take a structure $\widehat{J}$ of $\Im _{1}^{1}(T\overline{M})\;$defined by

\begin{equation}
\widehat{J}={\mathbf{F}}^{c}-\sum_{\alpha =1}^{r}(\xi _{\alpha }^{v}\otimes
\eta ^{\alpha v}-\xi _{\alpha }^{c}\otimes \eta ^{\alpha c})  \label{2.14}
\end{equation}

Similarly proof of \textbf{Theorem 4.1,} using (\ref{2.2}), (\ref{2.3}), (%
\ref{2.4}) and (\ref{2.5})and (\ref{2.6}) we have the equation

\begin{equation}
(\widehat{J})^{2}=\epsilon I.  \label{2.15}
\end{equation}
Thus, for $\widehat{J}$ is an almost complex (or paracomplex) structure in $%
T(\bar{M})$, the proof is completed.

Considering lift properties of tensor fields and the equation (\ref{2.14}),
we

\[
\widehat{J}X^{v}=({\mathbf{F}}X)^{v}+(\eta ^{\alpha }(X))^{v}\xi _{\alpha
}^{c}
\]

\[
\widehat{J}X^{c}=({\mathbf{F}}X)^{c}-(\eta ^{\alpha }(X))^{v}\xi _{\alpha
}^{v}+(\eta ^{\alpha }(X))^{c}\xi _{\alpha }^{c}.
\]

where $X\in \chi (\overline{M}),X^{v},X^{c}\in \chi (T\overline{M}).$

For example, if $\eta ^{\alpha }(X)=0,$ we obtain

\[
\widehat{J}X^{v}=({\mathbf{F}}X)^{v},\,\widehat{J}X^{c}=({\mathbf{F}}X)^{c}.
\]

In view of the structure $({\mathbf{F}},\xi _{\alpha },\eta ^{\alpha })\,$%
and the equations (\ref{2.11}),(\ref{2.13}) and (\ref{2.14}) it follows
\[
\widehat{J}\xi _{\alpha }^{v}=\delta _{\beta }^{\alpha }\xi _{\alpha
}^{c}=\xi _{\beta }^{c}\,,\,\,\,\,\widehat{J}\xi _{\alpha }^{c}=\delta
_{\beta }^{\alpha }\xi _{\alpha }^{v}=\xi _{\beta }^{v}\,,\,\,\alpha ,\beta
=1,2,...,r.
\]

\subsection{Horizontal Lifts}

\textbf{Theorem 4.3.} Let $({\mathbf{F}},\xi _{\alpha },\eta ^{\alpha })$ be
an almost $r-$paracontact( or $r-$paracontact) structure in $\overline{M}$
with an affine connection $\nabla $. Then an almost complex (or paracomplex)
in $T(\overline{M})$ is given by
\[
\widetilde{J}^{*}={\mathbf{F}}^{h}+\sum_{\alpha =1}^{r}\xi _{\alpha
}^{v}\otimes \eta ^{\alpha v}-\xi _{\alpha }^{h}\otimes \eta ^{\alpha h}).
\]

\textbf{Proof:} Taking into consideration the equation given by (\ref{1.9})
and the horizontal lifts of complex (or paracomplex) tensor fields, we have
\begin{equation}
({\mathbf{F}}^{2})^{h}=(\epsilon I+\sum_{\alpha =1}^{r}\xi _{\alpha }\otimes
\eta ^{\alpha })^{h},  \label{2.16}
\end{equation}
\begin{equation}
({\mathbf{F}}^{h})^{2}=\epsilon I+\sum_{\alpha =1}^{r}\xi _{\alpha
}^{h}\otimes \eta ^{\alpha v}+\xi _{\alpha }^{v}\otimes \eta ^{\alpha h},
\label{2.17}
\end{equation}
and
\begin{equation}
{\mathbf{F}}^{h}(\xi _{\alpha }^{h})=0,{\mathbf{F}}^{h}(\xi _{\alpha
}^{v})=0,  \label{2.18}
\end{equation}
\begin{equation}
\eta ^{\alpha h}\circ {\mathbf{F}}^{h}=0,\eta ^{\alpha v}\circ {\mathbf{F}}%
^{h}=0.  \label{2.19}
\end{equation}
\begin{equation}
\eta ^{\alpha h}(\xi _{\beta }^{h})=0,\eta ^{\alpha h}(\xi _{\beta
}^{v})=\delta _{\beta }^{\alpha },\eta ^{\alpha v}(\xi _{\beta }^{h})=\delta
_{\beta }^{\alpha }.\,  \label{2.20}
\end{equation}
Given an element $\widetilde{J}$ of $\Im _{1}^{1}(T\overline{M})\;$defined
by
\begin{equation}
\widetilde{J}^{*}={\mathbf{F}}^{h}+{\sum_{\alpha =1}^{r}}(\xi _{\alpha
}^{v}\otimes \eta ^{\alpha v}-\xi _{\alpha }^{h}\otimes \eta ^{\alpha h}).
\label{2.21}
\end{equation}
Thinking of the above equations, it is clear that
\begin{equation}
(\widetilde{J}^{*})^{2}=\epsilon I.  \label{2.22}
\end{equation}
Finally, $\widetilde{J}^{*}$ is an almost complex (or paracomplex) structure
in $T(\overline{M})$. Thus the theorem is proved.

Taking care of lift properties of tensor fields and the equation (\ref{2.21}%
), we

\[
\widetilde{J}^{*}X^{v}=({\mathbf{F}}X)^{v}-(\eta ^{\alpha }(X))^{v}\xi
_{\alpha }^{h}
\]

\[
\widetilde{J}^{*}X^{h}=({\mathbf{F}}X)^{h}+(\eta ^{\alpha }(X))^{v}\xi
_{\alpha }^{v}-(\eta ^{\alpha }(X))^{h}\xi _{\alpha }^{h}.
\]

where $X\in \chi (\overline{M}),X^{v},X^{h}\in \chi (T\overline{M}).$

For example, if $\eta ^{\alpha }(X)=0,$ we obtain

\[
\widetilde{J}^{*}X^{v}=({\mathbf{F}}X)^{v},\widetilde{J}^{*}X^{h}=({\mathbf{F%
}}X)^{h}.
\]

In view of the structure $({\mathbf{F}},\xi _{\alpha },\eta ^{\alpha })\,$%
and the equations (\ref{2.18}),(\ref{2.20}) and (\ref{2.22}) it follows
\[
\widetilde{J}^{*}\xi _{\alpha }^{v}=\delta _{\beta }^{\alpha }\xi _{\alpha
}^{h}=\xi _{\beta }^{h}\,,\,\,\,\widetilde{J}^{*}\xi _{\alpha }^{h}=\delta
_{\beta }^{\alpha }\xi _{\alpha }^{v}=\xi _{\beta }^{v}\,,\,\,\alpha ,\beta
=1,2,...,r.
\]

\thinspace \thinspace \thinspace \thinspace \thinspace \thinspace

\textbf{Theorem 4.4.} Let $({\mathbf{F}},\xi _{\alpha },\eta ^{\alpha })$ be
an almost Lorentzian $r-$paracontact ( or $r-$paracontact) structure in $%
\overline{M}$ with an affine connection $\nabla $. Then structure
\[
\widehat{J}^{*}={\mathbf{F}}^{h}-\sum_{\alpha =1}^{r}\xi _{\alpha
}^{v}\otimes \eta ^{\alpha v}-\xi _{\alpha }^{h}\otimes \eta ^{\alpha h}).
\]

is an almost complex (or paracomplex) in $T(\overline{M}).$

\textbf{Proof:} It can easily shown to taking care of the proofs of the
above.

\section{\textbf{Corollary}}

Taking into consideration the above theorems, we conclude that when we
consider an almost Lorentzian $r-$contact (or $r$- paracontact) structure on
the base manifold, the structure defined on the tangent bundle $T(\overline{M%
})$ is an almost complex (or paracomplex).

\thinspace \thinspace \thinspace \thinspace \thinspace \thinspace

\thinspace \textbf{REFERENCES}

[1] Tekkoyun, M., Civelek, \c{S}.: On lifts of structures on
complex manifolds, Differential Geometry-Dynamics Systems Vol.5,
No.1, 59-64 (2003).

[2] Tekkoyun, M.: On horizontal lifts of complex structures,
Hadronic Journal Supplement (HJS), Vol.18, No.4, 411-424 (2003).

[3] Tekkoyun, M.: On Para-Euler Lagrange and para- Hamiltonian
equations, Physics Letters A, Vol.340, 7-11 (2005).

[4] Tekkoyun, M.: On lifts of paracomplex structures, Turk. J.
Math., Vol.30, 197-210 (2006).

[5] Tekkoyun, M., Gorgulu,A.: Higher order complex Lagrangian and
Hamiltonian mechanics systems, Phyics Letters A, Vol.357, 261-269
(2006).

[6] Yano, K., Ishihara, S.: Almost complex structures induced in
tangent bundles, Kodai Math. Sem. Rep., Vol.19, 1-27(1967).

[7] Yano, K., Kobayashi, S.: Prolongations of tensor fields and
connections to tangent bundles. II. Infinitesimal automorphisms,
J. Math. Soc. Japan, Vol.18, 236-246(1966).

[8] Yano, K., Kobayashi, S.: Prolongations of tensor fields and
connections to tangent bundles. III. Holonomy groups, J. Math.
Soc. Japan, Vol.19, 486-488 (1967).

[9] Das, L.S., Khan, M.N.I.: Almost r-contact structures on the tangent
bundle, Differential Geometry-Dynamics Systems Vol.7, pp.34-41(2005).

[10] Blair, D.E.: Contact manifolds in Riemannian geometry,
Lecture Notes in Math, 509, Springer Verlag, New York, (1976).

[11] G\^{i}r\c{t}u, M.: An almost paracontact structure on the
indicatrix bundle of a Finsler space, Balkan Journal of Geometry
and Its Applications, Vol.7, No.2, 43-48 (2002).

[12] Ledger, A.J., Yano, K.: Almost complex structures on complex
bundles, Jour. Diff. Geometry, Vol.1, 355-368(1967).

[13] Matsumoto, K.: On Lorentzian paracontact manifolds, Bull. of
Yamagata Univ. Nat. Sci, Vol.12, No.2, 151-156 (1989).

[14] Niki\'{c}, J.: Conditions for invariant submanifold of a
manifold with the $(\varphi ,\xi ,\eta ,G)-$structure, Kragujevac
J. Math., Vol.25, 147-154 (2003).\\

Mehmet TEKKOYUN\\
Department of Mathematics, Faculty of Science and
Art,\\
Pamukkale University,\\
 20070 Denizli-TURKEY\\
e-mail:tekkoyun@pau.edu.tr
\end{titlepage}
\end{document}